\documentclass[]{IEEEtran}
\IEEEoverridecommandlockouts

\usepackage[utf8]{inputenc}
\usepackage[T1]{fontenc}
\usepackage[numbers]{natbib}
\usepackage{url}
\usepackage{physics} 
\usepackage{graphicx}
\usepackage{amsmath}
\usepackage[ruled,vlined]{algorithm2e}
\usepackage{amssymb,setspace,amsthm,calrsfs}
\usepackage{pgfplotstable,pstool,psfrag}
\pgfplotsset{compat=1.15}
\usepackage{float}
\usepackage{tabularx}
\usepackage{ltablex}
\usepackage{booktabs}
\usepackage{picins}
\usepackage{tabularx}
\usepackage{placeins}
\usepackage{supertabular}
\usepackage{tikz-cd}
\usepackage{subcaption}
\usepackage{xcolor}
\usepackage{color}
\definecolor{myblue}{RGB}{0,128,128}
\usepackage{acronym}

\newcommand{\R}{\mathbb{R}}
\newcommand{\N}{\mathbb{N}}

\newcommand{\probability}{\mathbb{P}}

\def\defeq{\mathrel{\mathop:}=}
\newtheorem{thm}{Theorem}
\newtheorem{dfn}{Definition}

\newtheorem{rem}{Remark}

\makeatletter
\let\mytagform@=\tagform@
\def\tagform@#1{\maketag@@@{\color{myblue}(#1)}}
\newif\ifnoindentafter
\noindentafterfalse
\makeatother

\newif\iffinalfonts
\iffinalfonts
\else

\fi

\SetAlFnt{\small}
\SetAlCapFnt{\small}

\usepackage[urlcolor=myblue, citecolor=myblue, colorlinks=true,linkcolor=myblue,urlcolor=myblue,citecolor=myblue]{hyperref}
\hypersetup{linkbordercolor=myblue}
\SetAlCapNameFnt{\small}
\SetAlCapHSkip{0pt}
\IncMargin{-\parindent}
\makeatletter
\newenvironment{smalleralign}[1][\small]
 {\par\nopagebreak\leavevmode\vspace*{-\baselineskip}%
  \skip0=\abovedisplayskip
  #1%
  \def\maketag@@@##1{\hbox{\m@th\normalfont\normalsize##1}}%
  \abovedisplayskip=\skip0
  \align}
 {\endalign\ignorespacesafterend}
\makeatother

\begin{document}

\title{Predictability and Fairness in Load Aggregation with Deadband}



\author{
F. V. Difonzo,
M. Roubal{\' i}k,
J. Marecek
\thanks{J. Marecek is at CEZ a.s., Prague, the Czech Republic.}
\thanks{M. Roubalik is at StormGeo-Nena, 
Oslo, Norway.}
\thanks{J. Marecek and M. Roubalik are at the Czech Technical University, Prague, the Czech Republic.}
\thanks{F. V. Difonzo is at the Instituto per le Applicazioni del Calcolo "Mauro Picone" (IAC-CNR), Italy.}
\thanks{Manuscript received \today.}
}
\maketitle

\begin{abstract}
Virtual power plants and load aggregation are becoming increasingly common. There, one regulates the aggregate power output of an ensemble of distributed energy resources (DERs).
Marecek et al. [Automatica, Volume 147, January 2023, 110743] recently suggested that long-term averages of prices or incentives offered should exist and be independent of the initial states of the operators of the DER, the aggregator, and the power grid.
This can be seen as predictability, which underlies fairness.
Unfortunately, the existence of such averages cannot be guaranteed with many traditional regulators, including the proportional-integral (PI) regulator with or without deadband. 
Here, we consider the effects of losses in the alternating current model and the deadband in the controller.
This yields a non-linear dynamical system (due to the non-linear losses) exhibiting discontinuities (due to the deadband).
We show that Filippov invariant measures enable reasoning about predictability and fairness while considering non-linearity of the alternating-current model and deadband. 
\end{abstract}


\section{Introduction}
Load aggregation of distributed energy resources (DERs) is an increasingly common response to the increasing integration of intermittent sources of energy \cite{lund2015review}.
In the European Union, this is now mandated by Regulation (EU) 2019/943 on the internal market for electricity sets, specifically requiring ``fair rules''.
More recent plans known as the RePowerEUR (COM/2022/230) suggest ``fairness and solidarity are defining principles''.
The fairness \cite{marecek2021predictability} may suggest that similar operators of DER receive a similar treatment in terms of load reductions, disconnections, incentives, prices, or power quality, over the long term, independent on the initial conditions, such as consumption at some point in the past. Below, we define these criteria formally in terms of properties of a certain stochastic model. 

\begin{figure}[t!]
\includegraphics[width=\columnwidth]{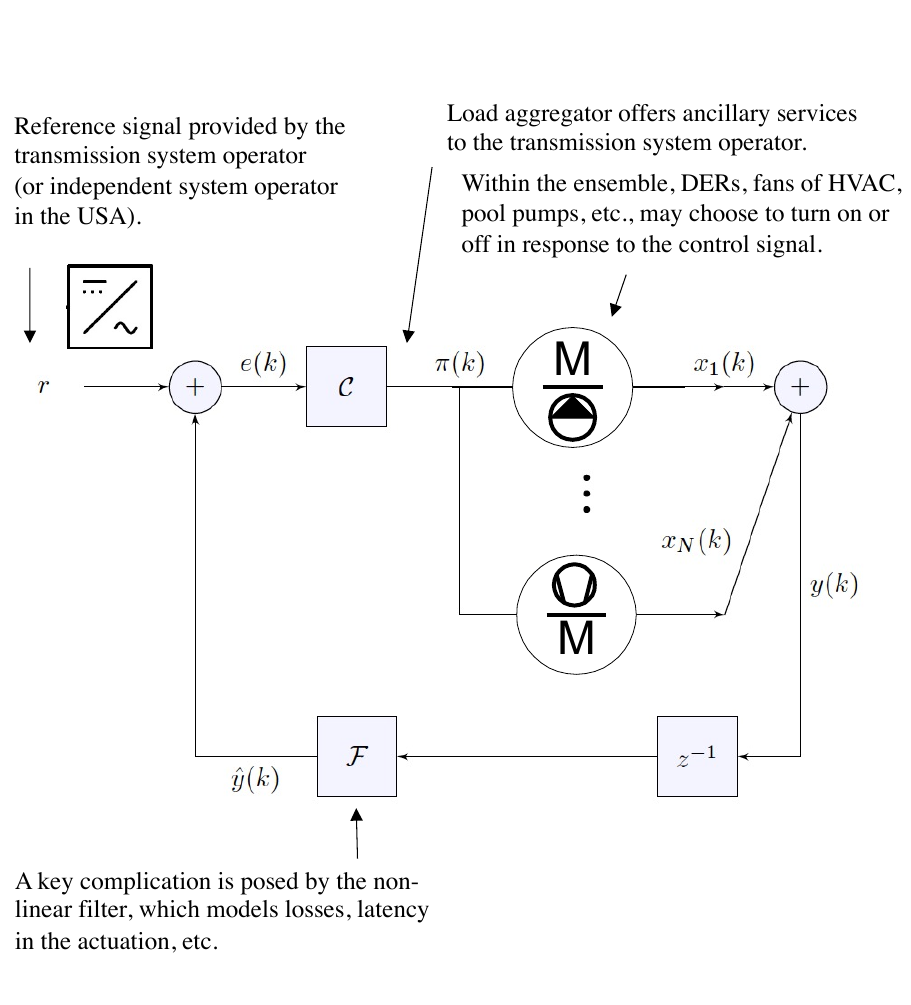}
\caption{An illustration of a closed-loop model, taken from \cite{marecek2021predictability}. In this paper, we extend the reasoning to discontinuous dynamical systems obtained by considering deadband in the controller $\mathcal C$.}
\label{system}
\end{figure}

The regulation of the aggregate power output \citep[e.g.]{Hiskens2011,pinson2014benefits} is challenging for a number of reasons.
First, it needs care in modelling participants of various forms, including roof-top photovoltaic systems, interruptible or deferrable loads, and battery energy storage systems.
Second, as highlighted by the 2030 roadmap of the IEEE Control Systems Society \cite[cf. Section 2.D]{alleyne2023control}, it also presents an important technical challenge in control theory. 
The load aggregation introduces a new feedback loop \cite{6197252,li2020real} into a system. 
The modelling of the non-deterministic behaviour of a number of DERs, which may be time-varying,
 is challenging particularly when one wishes to consider ``microscopic'' properties,
 at the level of one DER. 
This includes ``predictability and fairness'' properties of Fioravanti et al. \cite{ErgodicControlAutomatica,marecek2021predictability},
which require certain ergodicity. (Cf. page 43 in \cite{alleyne2023control}.) 
Besides, the system is already non-linear due to the losses in AC transmission 
and the use of deadband renders the dynamics discontinuous. 
Recent studies of load aggregation by Li et al. \cite{li2020real}, Marecek et al. \cite{marecek2021predictability},
and Kasis et al. \cite{KASIS2020104755,Kasis2021} have made a considerable progress in addressing the the closed the closed-loop aspects and the non-linearity of the losses.

In this paper, we aim to address all the challenges concurrently, including the deadband. 
Mathematically, this requires a fundamentally different approach, utilizing stochastic difference inclusions, 
rather than stochastic difference equations \cite[as in]{ErgodicControlAutomatica,marecek2021predictability} or difference inclusions \cite[as in]{KASIS2020104755,Kasis2021}.
We present a framework for reasoning about ergodic properties of closed-loops with non-smooth controllers, such as deadband, and non-smooth response of the agents. 
We illustrate the ergodic behaviour in simulations using Matpower \cite{zimmerman2010matpower} on standard IEEE test cases and a test case of \cite{marecek2021predictability}.

Practically, this is important for two reasons. 
While deadband has been traditionally used primarily to protect the turbines of thermal and hydropower
generators, the DERs such as pool pumps and compressors of HVAC are likewise quite sensitive to short-cycling. 
Enforcing a no-short-cycling requirement \cite{7345604} may resemble a deadband in all but name.  

Second, predictability and fairness cannot be guaranteed with many traditional controllers used by the load aggregator, including the usual proportional-integral (PI) controller, with or without a deadband. 
In particular, in the model of Fioravanti et al. \cite{ErgodicControlAutomatica,alleyne2023control}, the closed-loop system cannot have a unique invariant measure, making it it impossible to satisfy the fairness properties, 
for any controller ${\mathcal C}$ that is linear and marginally stable with a pole $s_1 =e^{q j \pi}$ on the unit circle, where $q$ is a rational number. 
With deadband, one introduces a discontinuity into the dynamics, which complicates the analysis substantially,
but preserves this negative result. 

\section{Definitions and Related Work}\label{sec:model}



\subsection{Related Models}\label{ssec:model}

As in \cite{ErgodicControlAutomatica,marecek2021predictability}, we consider a discrete-time, synchronous model of an ensemble of DERs illustrated in Figure \ref{system}. We model these as discrete-time stochastic inclusions: 

\begin{dfn}[Discrete-time inclusion, \cite{smirnov}]
Let $F$ be a set-valued map, whose values are closed subsets of $\mathbb{R}^n$. Then, \emph{a discrete-time inclusion} is defined by:
\begin{align}\label{eq:dfi}
    x(k+1)\in F(x(k)),\quad k\geq0,
\end{align}
for some sequence $\{x(k)\}_{k\in\mathbb{N}}\subseteq\mathbb{R}^n$.
\end{dfn}

Notice that Markov chain can be seen as a natural extension, which assigns positive probability of moving from $x(k)$ to states within $F(x(k))$.
A key notion across both Markov chains and other stochastic extensions of difference inclusions is the invariant measure  $\mathcal P \mu = \mu$ for some transition operator $\mathcal P$.
In Section \ref{sec:nonlinear}, we generalize these notions to account for set-valued maps $F$ in \eqref{eq:dfi} and their Filippov \cite{filippov} convexifications.

Returning to the model of \cite{ErgodicControlAutomatica,marecek2021predictability}, the same signal $\pi(k)\in \Pi \subseteq \mathbb{R}$ is sent synchronously to all $N$ DERs at time $k$. Let us now adapt the model to the discontinuous setting, necessary for what follows, and consider two $C^2$-functions $h_W,h_H:\R^n\to\R$, both having $0$ as regular value.
Each DER $i$ has a state $x_i\in \R^{n_i}$ and we assume that there are $w_i \in \N$ state transition maps, discontinuous at $h_W=0$,
\begin{equation}\label{eq:W}
\mathcal{W}_{ij}(x)=
\begin{cases}
\mathcal{W}_{ij}^-(x), & h_W(x)<0, \\
\mathcal{W}_{ij}^+(x), & h_W(x)>0,
\end{cases}\quad
j=1,\ldots,w_i,
\end{equation}
for agent $i$, and $h_i \in \N$ output maps, discontinuous at $h_H=0$,
\begin{equation}\label{eq:H}
\mathcal{H}_{il}(x)=
\begin{cases}
\mathcal{H}_{il}^-(x), & h_H(x)<0, \\
\mathcal{H}_{il}^+(x), & h_H(x)>0.
\end{cases}\quad
l=1,\ldots,h_i,
\end{equation}
for the same agent $i$.
One can model the system using a stochastic difference inclusion:
\begin{align}\label{eq:nonlinear-agents_DiffIncl}
x_i(k+1) & \in  \{  {\mathcal W}_{ij}^\pm(x_i(k)) \;\vert\; j = 1, \ldots, w_i\}, \\
y_i(k) & \in  \{ {\mathcal H}_{i\ell}^\pm(x_i(k)) \;\vert\; \ell = 1, \ldots, h_i\},
\end{align}
governed by probability functions $p_{ij}^\pm : \Pi \to [0,1]$, $j=1,\ldots,w_i$, and
$p'^\pm_{i\ell} : \Pi \to [0,1]$, $\ell=1,\ldots,h_i$, respectively, where for all $\pi \in \Pi$, $i=1,\ldots,N$ 
$\sum_{j=1}^{w_i} p_{ij}^\pm(\pi) = \sum_{\ell=1}^{h_i} p'^\pm_{i\ell}(\pi) = 1$.
This means that:
\begin{subequations} \label{eq:problaws}
\begin{align}\label{eq:problaw-1}
&\mathbb{P}\big( x_i(k+1) = {\mathcal W}_{ij}^\pm(x_i(k)) \big) = p_{ij}^\pm(\pi(k)),\\
&\mathbb{P}\big( y_i(k) = {\mathcal H}_{i\ell}^\pm(x_i(k)) \big) = p'^\pm_{i\ell}(\pi(k))
\end{align}
\end{subequations}
for all agents $i$ and time steps $k\in\N$. 
The outputs $y_i(k)$ are Markovian, in that they depend only on the most recent states $x_i(k)$ and the most recent signal $\pi(k)$.

To complete a Markovian model, one should also consider the aggregate output $y(k) = \sum_{i=1}^N y_i(k)$, the state $x_f(k)$ of a non-linear filter:
\begin{align}
\label{eq:nonlinear-filter}
&{\mathcal F} ~:~ \begin{cases}
x_f(k+1) = {\mathcal W}_{f}(x_f(k),y(k)), & \\
\hat y(k) = {\mathcal H}_{f}(x_f(k),y(k)), &
\end{cases}
\end{align}
and the state $x_c(k)$ of the non-linear controller: 
\begin{align}
\label{eq:nonlinear-cont}
&{\mathcal C} ~:~ \begin{cases}
x_c(k+1) = {\mathcal W}_{c}(x_c(k),\hat y(k),r), & \\
\pi(k) = {\mathcal H}_{c}(x_c(k),\hat y(k),r). &
\end{cases}
\end{align}
Let us use $\mathbb{X}_i, i=1,\ldots,N, \mathbb{X}_C$ and $\mathbb{X}_F$ to denote the state spaces of the
agents, the controller and the filter, respectively.
Then, one can consider the evolution of a system on state space $\mathbb{X} \defeq \prod_{i=1}^N \mathbb{X}_i \times \mathbb{X}_C \times \mathbb{X}_F$
governed by: 
\begin{equation}
\label{eq:6}
x(k+1) \defeq \begin{pmatrix}
(x_i)_{i=1}^N  \\
x_f \\
x_c
\end{pmatrix} (k+1) \in \{ F_m^\pm(x(k)) \,\vert\, m \in {\mathbb M}\},
\end{equation}
where the maps $F_m^\pm: \mathbb{X} \to \mathbb{X}$ are of the form
\begin{equation}
\label{eq:F_m-definition}
F_m^\pm(x(k)) \defeq  \begin{pmatrix}
( {\mathcal W}_{ij}^\pm(x_i (k)) )_{i=1}^N  \\
{\mathcal W}_f(x_f(k), \sum_{i=1}^N  {\mathcal H}_{i\ell}^\pm(x_i (k))) \\
{\mathcal W}_c(x_c(k), {\mathcal H}_f(x_f(k), \sum_{i=1}^N  {\mathcal H}_{i\ell}^\pm(x_i (k))))
\end{pmatrix}.
\end{equation}
The index sets for $F_m^\pm$ coincide and are given by: 
\begin{equation}
\label{eq:8}
\mathbb{M} \defeq \prod_{i=1}^N \{ (i,1), \ldots, (i,w_i) \} \times \prod_{i=1}^N \{ (i,1), \ldots, (i,h_i) \}
\end{equation}
where for each multi-index
$m=((1,j_1),\ldots,(N,j_N),(1,l_1),\ldots,(N,l_N))$, the
probability of choosing $F_m^\pm$ is:
\begin{multline}
\label{eq:9}
\probability \left( x(k+1)=F_m^\pm(x(k)) \right) = \\ \left(\prod_{i=1}^N
p_{ij_i}^\pm(\pi(k)) \right) \left( \prod_{i=1}^N p'^\pm_{il_i}(\pi(k))
\right) =: q_m^\pm(\pi(k)).
\end{multline}

\subsection{Related Results}


Fioravanti et al. \cite{ErgodicControlAutomatica} considered the regulation of a generalization of load aggregation
in the feedback system of Figure \ref{system}, but again without the non-linear filters to model alternating-current transmission losses. 
They consider predictability, in the sense that for each agent $i$ there exists a
constant $\overline{r}_i$ such that
\begin{equation}
\label{eq:3}
\lim_{k\to \infty} \frac{1}{k+1} \sum_{j=0}^k y_i(j) = \overline{r}_i,
\end{equation}
where the limit is independent of initial conditions, and of fairness in the sense that  the limits $\overline{r}_i$ coincide
for all $1 \le i \le N$.
They showed that if the transfer function of the controller ${\mathcal C}$ is a single-input single-output (SISO) system with a pole $s_1 =
e^{q j \pi}$ on the unit circle where $q$ is a rational number, i.e., is not asymptotically stable, the closed-loop system cannot be uniquely ergodic, and thus it can be neither predictable nor fair. Thus, even though the PI controller (or any other controller an integral action) performs its regulation function well, 
it destroys the ergodic properties of the closed loop.

In a follow-up study, Marecek et al. \cite{marecek2021predictability} extended the results to the non-linearity of the alternating-current model and considered the ergodic properties of incremental input-to-state controllers and filters. In doing so, the authors utilize {\em iterated function systems} \citep{elton1987ergodic,BarnsleyDemkoEltonEtAl1988,barnsley1989recurrent} following \cite{ErgodicControlAutomatica}. Their main result can be rephrased as follows.
In the feedback system of Figure \ref{system},
the state of agent 
$i \in \{1,\cdots,N\}$ is governed by
\begin{subequations}\label{eq:nonlinear-agents}
\begin{align}
x_i(k+1) &= {\mathcal W}_{ij}(x_i(k)) \label{eq:nonlinear-agentsW}\\
y_i(k) &= {\mathcal H}_{ij}(x_i(k)), \label{eq:nonlinear-agentsH}
\end{align}
\end{subequations}
where functions ${\mathcal W}_{ij}$ and ${\mathcal H}_{ij}$ are continuously differentiable and globally Lipschitz-continuous with Lipschitz constants $l_{ij}$, and $l'_{ij}$, respectively. 
Probability functions $p_{ij},p'_{il}$ are  Dini continuous and bounded from blow by scalars $\delta, \delta' > 0$ for all $\pi\in \Pi$ and all $(i,j)$.
When agents' transition maps and the probability functions satisfy a certain contraction-on-average condition, then for every incrementally input-to-state stable controller $\mathcal{C}$ and every incrementally input-to-state stable filter $\mathcal{F}$ compatible with the feedback structure, the feedback loop has a unique invariant measure. 

There, Marecek et al. \cite{marecek2021predictability} use one of the standard definitions of incremental input-to-state stability (ISS), which we reiterate here for the convenience of the reader:

\begin{dfn}[Comparison functions] \label{df:class-K}
A function $\gamma : \mathbb R^{+}\to \mathbb R^{+}$ is is said to be of class $\mathcal K$ if it is continuous, increasing and $\gamma(0)=0$. It is of class $\mathcal K_{\infty}$ if, in addition, it is proper, i.e., unbounded.
A continuous function $\beta: \mathbb R^{+}\times \mathbb R^{+}\to \mathbb R^{+}$ is said to be of class $\mathcal K \mathcal L$, if for all fixed $t$ the function $\beta(\cdot, t)$ is of class $\mathcal K$ and for all fixed $s$, the function $\beta(s, \cdot)$ is non-increasing
and tends to zero as $t\to \infty$.
\end{dfn}

\begin{dfn}[Incremental ISS, \cite{angeli2002lyapunov}]
Let $\mathcal U
$ denote the set of all 
input functions $u: \mathbb Z_{\ge k_0}\to \mathbb R^d$%
. Suppose $F: \mathbb R^d \times \R^n\to \mathbb R^n$ is continuous, then the discrete-time non-linear dynamical system
\begin{align}
x(k+1) = F(x(k),u(k))
\label{system-inciss}
\end{align}
is called (globally) \emph{incrementally input-to-state-stable} (incrementally ISS), if there exist $\beta\in \mathcal{KL}$ and $\gamma\in\mathcal K$ such that
for any pair of inputs $u_1, u_1\in\mathcal{U}$ and any pair of initial condition  $\xi_1, \xi_2 \in \R^n$:
\begin{align*}
&\|x(k, \xi_1, u_1)-x(k, \xi_2, u_2)\|\\ 
&\le \beta(\|\xi_1 - \xi_2\|, k)+ \gamma(\|u_1-u_2\|_{\infty}),\quad \forall k\in \N.   
\end{align*}
\end{dfn}

We refer to \cite{angeli2002lyapunov} for a beautiful exposition.

\begin{figure}
\centering
\begin{subfigure}[b]{0.4\textwidth}
    \caption{Evolutions of active powers $(P_1, P_2)$ over time $t$.}
    \begin{center}
    \includegraphics[width=0.7\columnwidth]{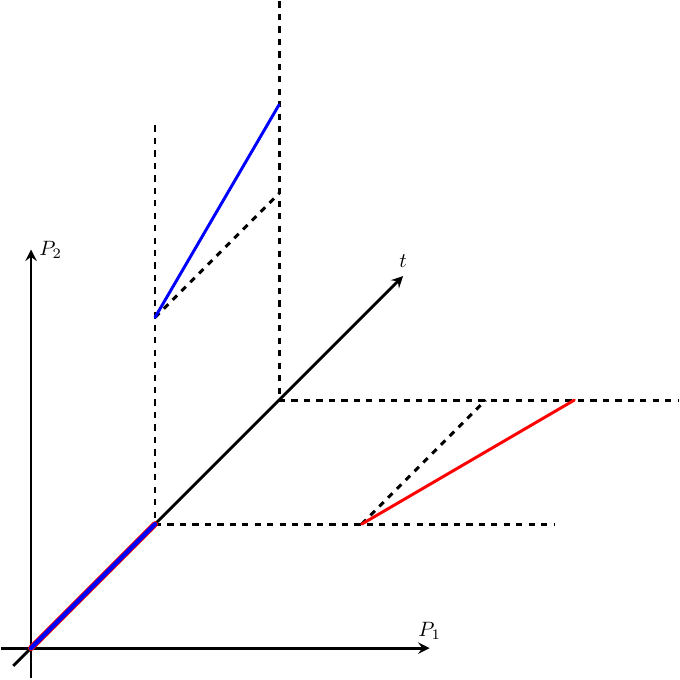}
    \end{center}
    \label{fig:filippov_a}
\end{subfigure}
\\[5pt]
\begin{subfigure}[b]{0.35\textwidth}
    \caption{Neighborhoods of evolutions of active powers $(P_1, P_2)$ over time $t$.}
    \begin{center}    
    \includegraphics[width=0.7\columnwidth]{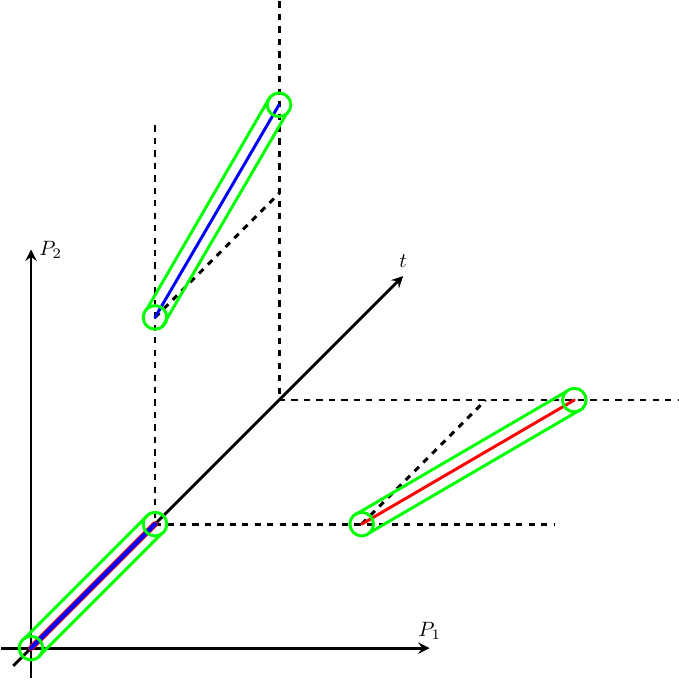}
    \end{center}
    \label{fig:filippov_b}
\end{subfigure} 
\\[5pt]
\begin{subfigure}[b]{0.4\textwidth}
    \caption{Convex hulls of the neighborhoods of evolutions of active powers $(P_1, P_2)$ over time $t$.}
    \begin{center}    
    \includegraphics[width=0.7\columnwidth]{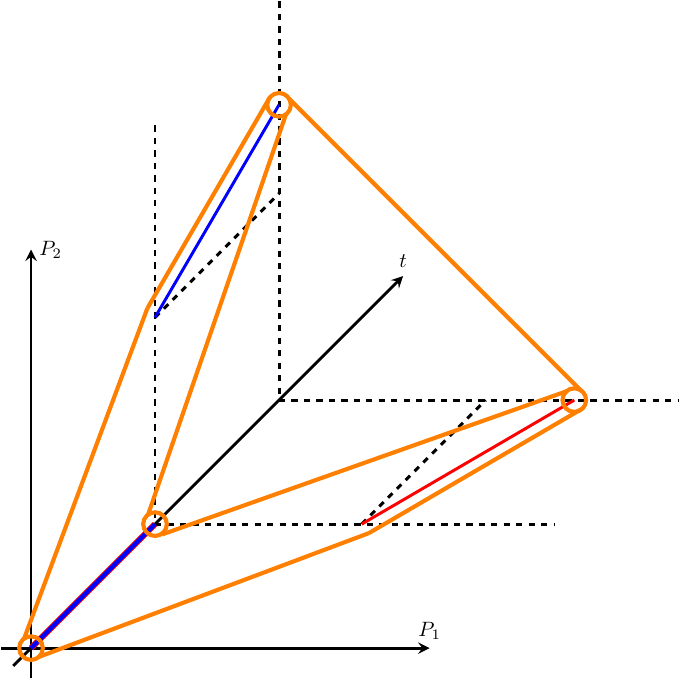}
    \end{center}        
    \label{fig:filippov_c}
\end{subfigure}
\caption{An illustration of the Filippov convexification in three easy steps. On the top of the column, (a) displays two trajectories of discontinuous dynamics of active power output over time, given by deadband.
Below, (b) constructs neighbourhoods around the continuous subsets of the same trajectories. At the bottom, (c) suggests the convexification.}
\label{fig:filippov}
\end{figure}

\section{Reasoning about Predictability and Fairness with Deadband}\label{sec:nonlinear}

It is natural to investigate how to generalize the assumptions of Marecek et al. \cite{marecek2021predictability} to accommodate controllers with a deadband or other discontinuitites. 
In particular, it could be that output 
depends on the state through some piecewise smooth control function $\phi$. For example, the control function $\phi$ can represent a deadband: 
\begin{align}
y_i(k)=\mathcal{H}_{ij}(\phi(x_i(k))),
\end{align}
where
\begin{align}
\label{eq:deadband}
\phi(x):=
\begin{cases}
x, & |x|\geq\delta, \\
0, & |x|<\delta,
\end{cases}
\end{align}
for some $\delta>0$. This would imply that assumption $(i)$ in Theorem 2 of \cite{marecek2021predictability} should be weakened in order to retain, in some extended sense, an attractive invariant measure and unique ergodicity for the resulting system. More specifically, we have to require $\mathcal{W}_{ij},\mathcal{H}_{ij}$ to be Lebesgue measurable and $L^{\infty}_{\textrm loc}$ (see \cite{Haddad2014}). In doing so, 
we lose smoothness. 
Thus, classical solutions need not exist, and one  has to extend the concept of solution.

%

A common approach is to consider the operator $F[\cdot]$ in \eqref{eq:nonlinear-agents_DiffIncl} to be the Filippov convexification of the original flow:
\begin{dfn}[Filippov convexification, \cite{filippov}] \label{def:FilippovConv}
Let $f:\R^{n}\to\R^{n}$ be Lebesgue measurable and $L^{\infty}_{\textrm loc}$, and let $x\in\R^{n}$. The \emph{Filippov convexification} of $f$ at $x$ is defined as the set-valued map 
\begin{equation}\label{eq:FilippovF}
F[f](x):=\bigcap_{\delta>0}\bigcap_{\mu(S)=0}\overline{\textrm{co}}\{f(B(x,\delta))\setminus S\},
\end{equation}
where $\mu$ is the Lebesgue measure on $\R^{n}$, $\overline{\textrm{co}}(A)$ is the closed convex hull of a generic set $A$ and $B(x,\delta)$ is the open ball centered at $x$ with radius $\delta$ (see \cite{filippov} for further details).
\end{dfn}

We refer to Figure \ref{fig:filippov} for an illustration of the Filippov convexification on a simple example of a 2D example, illustrating 
the evolution of active-power output of two DERs. 
In Figure \ref{fig:filippov_a}, there are plotted power outputs of two generators ($P_1$ and $P_2$) over time $t$. At first, both generators produce the same amount of active power, since they are both off. At certain moment, both generators are switched on. This leads to obvious discontinuity in both $P_1$ and $P_2$ as functions of time.
Figure \ref{fig:filippov_b} shows neighborhoods of continuous subsets of both $P_1$ and $P_2$ trajectories. Then, in Figure \ref{fig:filippov_c}, convex hull around the neighborhoods is constructed.

The set-valued map, which results from \eqref{eq:FilippovF}, is known to be upper semicontinuous with closed convex values, and Lipschitz continuous under specific assumptions (see \cite{smirnov}). See \cite{Cortes08discontinuousdynamical,CortesBullo2009} for a thorough exposition on the topic.
The first of these properties guarantees that \eqref{eq:nonlinear-agents_DiffIncl} admits at least a Filippov solution:

\begin{dfn}[Filippov solution, \cite{filippov}] \label{def:FilippovSol}
A \emph{Filippov solution} to \eqref{eq:nonlinear-agents_DiffIncl} is a couple of sequences $\{(x_i(k),y_i(k))\}_{k\in\mathbb{N}}$ satisfying, for each $k\in\mathbb{N}$,
\begin{subequations}\label{eq:nonlinear-disc}
\begin{align}
    x_i(k+1) & \in  F[{\mathcal W}_{ij}](x_i(k)),\,j=1,\ldots,w_i, \\
y_i(k) & \in  F[{\mathcal H}_{i\ell}](x_i(k)),\,l=1,\ldots,h_i,
\end{align}
\end{subequations}
with $F[\cdot]$ standing for the Filippov convexification of the corresponding flow.
\end{dfn}

Then, using the concepts of dissipativity for discontinuous dynamical systems with appropriate storage functions and supply rates (see \cite{HaddadSadikhov2014}) one can then try to extend results in Theorem 2 of \cite{marecek2021predictability}  to \eqref{eq:nonlinear-agents_DiffIncl}.
This way, one could specialise the notion of a Filippov invariant measure \cite{NOVAES2021102954}: 

\begin{dfn}[Filippov invariant measure] \label{def:FilippovInvariantMeasure}
We say that the discrete-time inclusion system \eqref{eq:dfi} preserves a measure $\mu$ on $\R^n$ if
\[
\mu(\mathcal{S}_F(A)(k))=\mu(A),\quad k\in\mathbb{N},
\]
for any Borel subset $A\subseteq\R^n$, where
\[
\mathcal{S}_F(A)(k):=\bigcup_{x\in A}\mathcal{S}_F(x)(k),
\]
where $\mathcal{S}_F(x)$ denotes the set of all maximal solutions to \eqref{eq:dfi}, for a given initial condition.
\end{dfn}

Now, within the above setting, one could then prove existence and uniqueness of invariant measures for the resulting Filippov system.
In turn, the existence of the unique invariant measure guarantees predictability. 

\begin{rem}
    Other types of convexification could be considered in \eqref{eq:nonlinear-agents_DiffIncl}, such as Krasovskii convexification, which could be useful for time-varying non-convex optimization problems \cite{Hauswirth2018TimevaryingPD}. However, possibly less is known about invariant measures for the corresponding resulting system than for the Filippov convexification, making them less attractive in the current context.
\end{rem}

As a successive analysis, following
\cite{HADDADHui2009,HaddadSadikhov2015}, optimality for time-invariant discontinuous control systems on the infinite interval utilizing a steady-state Hamilton-Jacobi-Bellman approach for characterizing optimal discontinuous nonlinear feedback controllers can also be considered, together with results gain, sector and disk margin, relative to some suitable performance criterion.

\section{The Main Result}

Our main result concerns the ergodic properties of non-linear controllers and filters in the discontinuous setting:

\begin{thm}
\label{thm:main}
Consider the feedback system depicted in Figure \ref{system}.
Assume that each agent
$i \in \{1,\cdots,N\}$ has a state governed by the non-linear iterated function system \eqref{eq:nonlinear-disc}
where, for some event $C^2$-functions $h_W,h_H:\R^n\to\R$, with $0$ as a regular value, 
both \eqref{eq:W} and \eqref{eq:H} hold. \\
Let $h:=h_W=h_H$ and $\Sigma^\pm:=\{x\in\R^n\,:\,h(x)\gtrless0\}$. Let us further assume that ${\mathcal W}_{ij}^\pm$, ${\mathcal H}_{il}^\pm$ and the compositions $F_m^\pm$ 
\eqref{eq:F_m-definition} 
are such that:
\begin{itemize}
\item[(i)] 
${\mathcal W}_{ij}^\pm$ and ${\mathcal H}_{il}^\pm$ are globally $l_{ij}$-Lipschitz and $l'_{il}$-Lipschitz-continuous, respectively, and $C^1$-functions;
\item[(ii)]
$p_{ij}^\pm,p'^\pm_{il}$ are Dini continuous probability functions so that
the probabilistic laws \eqref{eq:problaws} are satisfied;
\item[(iii)] 
there exists $\delta, \delta' > 0$ such that
$p_{ij}^\pm(\pi) \geq  \delta > 0$,
$p'^\pm_{il}(\pi) \geq  \delta' > 0$ for all $\pi\in \Pi$ and all $(i,j)$;
\item[(iv)] 
for each multi-index $m$, the composition $F_{a,m}^\pm: \prod_{i=1}^{N} \mathbb {X}_i \to \prod_{i=1}^{N} \mathbb {X}_i$,  $F_{a,m}^\pm\left((x_i)_{i=1}^N\right):=\left(\mathcal{W}_{im}^\pm(x_i)\right)_{i=1}^N$ 
of agents' transition maps and the probability functions $p_m^\pm$ satisfies the contraction-on-average condition, i.e., there exists a constant $0<\tilde c<1$ such that for all $x,\hat{x} \in \mathbb{X}_a, x\neq\hat{x}, \pi\in\Pi$ we have
\begin{equation}
\label{eq:average-contraction}
\sum_{m}p_m^\pm(\pi) \frac{\|F_{a,m}^\pm(x)-F_{a,m}^\pm(\hat{x})\|}{\|x-\hat{x}\|} < \tilde c.
\end{equation}
.
\end{itemize}
Moreover, let us assume that for every invariant measures  $\nu^\pm$ of the feedback loop relative to ${\mathcal W}_{ij}^\pm$, ${\mathcal H}_{il}^\pm$, there exist a measure $\mu$ and constant values $\alpha^\pm$ such that $\nu^\pm=\alpha^\pm\mu$ on $\Sigma^\pm$ and, for all $x\in\Sigma:=\mathrm{bnd}\Sigma^+\cap\mathrm{bnd}\Sigma^-$,
\begin{align*}
\alpha^+\nabla h(x)^\top\mathcal{W}_{ij}^+ &= \alpha^-\nabla h(x)^\top\mathcal{W}_{ij}^-, \\
\alpha^+\nabla h(x)^\top\mathcal{H}_{il}^+ &= \alpha^-\nabla h(x)^\top\mathcal{H}_{il}^-.
\end{align*}
Then, for every incrementally input-to-state stable controller $\mathcal{C}$ and every incrementally input-to-state stable filter $\mathcal{F}$ compatible with the feedback structure, the feedback loop has a unique Filippov invariant measure.
\end{thm}
\begin{proof}
Since $\mathcal{W}_{ij}^\pm,\mathcal{H}_{il}^\pm$ are locally essentially bounded, from classical results in \cite{filippov} it follows that both $F[\mathcal{W}_{ij}]$ and $F[\mathcal{H}_{il}]$ are upper semicontinuous and have nonempty, compact and convex values. Moreover, from the assumptions, for every incrementally input-to-state stable controller $\mathcal{C}$ and every incrementally input-to-state stable filter $\mathcal{F}$ compatible with the feedback structure, the feedback loop has unique Filippov invariant measures $\nu^\pm$. Thus there exist a measure $\mu$ and constant values $\alpha^\pm$ such that $\nu^\pm=\alpha^\pm\mu$. Then, from \cite[Theorem A]{NOVAES2021102954}, Filippov system \eqref{eq:nonlinear-disc} preserves the measure $\nu:=f\mu$, where
\[
f(x):=
\begin{cases}
\alpha^-, & h(x)<0, \\
\alpha^+, & h(x)>0.
\end{cases}
\]
If, by contradiction, there exists $\nu'\neq\nu$  Filippov invariant measure for \eqref{eq:nonlinear-disc}, then there would exist, from \cite[Theorem A]{NOVAES2021102954}, invariant measures ${\nu'}^\pm$ for the feedback loop relative to $\mathcal{W}_{ij}^\pm,\mathcal{H}_{il}^\pm$. By uniqueness, these measures must coincide with $\nu^\pm$, so that $\nu'=\nu$, which is a contradiction. The claim is proved.
\end{proof}

\begin{rem}
Theorem \ref{thm:main} implies that the piecewise dynamical system associated to \eqref{eq:nonlinear-disc} is of co-dimension $1$. This is needed in order to use results about invariant measures in Filippov system from \cite{NOVAES2021102954}. At the best of our knowledge, nothing is known about co-dimension $2$ Filippov systems, namely the case of general event smooth functions $h_W,h_H$. Further notions of ergodicity are an important question for further study, both across the power-systems applications and code non-smooth dynamical systems.
\end{rem}

\section{Simulating the Impact of Deadband}\label{sec:numeric}
Our numerical results are based on simulations utilising Matpower \cite{zimmerman2010matpower}, a standard open-source toolkit for power-systems analysis in version 7.1 running on Mathworks Matlab 2019b.
We use Matpower's power-flow routine to implement a non-linear filter, which models the losses in the alternating-current model. 

We consider two test cases:
\begin{itemize}
    \item First test case, closely related to \cite{marecek2021predictability}, with number of agents $N = 20$. The system is regulated to $r = N/2$. 
    \item Standard IEEE 30-bus system with assembly of $N = 50$ agents located at bus n. 2. We aim to regulate the system to the reference signal $r$, which is $60.97$ MW plus losses at time $k=0$.
\end{itemize}

The response of each system $\mathcal S_i$ to the control signal $\pi(k)$ provided by the load aggregator at time $k$ takes the form of probability functions, which represent the probability that DERs are committed at time $k$, of the control signal $\pi(k)$. The probability functions are defined as

\begin{align}\label{eq:prob-func}
g_{i1}(x_i (k+1)=1)&=0.02+\frac{0.95}{1+\exp(-\xi\pi(k)-x_{01})}
\end{align}
modelling the response of this first half of agents, whereas the second half of agents is modelled by probability function:
\begin{align}\label{eq:prob-func2}
g_{i2}(x_i (k+1)=1)&=0.98-\frac{0.95}{1+\exp(-\xi\pi(k)-x_{02})}
\end{align}

Initially, the first ($N/2$) generators are off and the second half of the generators are on.
The commitment of the DERs within the ensemble is provided to Matpower which, using a standard power-flow algorithm, computes 
the active power output $P(k)$ of the individual DERs:
\begin{align}\label{eq:output-gnrtr}
P(k)\in \mathbb{R}^N 
\end{align} 
taking into account losses in transmission. Thus computed active power output is aggregated into the total active power output $p(k) = \sum P(k)$ of the ensemble.
Then, a filter $\mathcal F$ is applied: 
\begin{align}\label{eq:filtr}
\hat p(k)= \frac{p(k)+p(k-1)}{2} + \textrm{losses}(k).
\end{align}
Notice that here we both smooth the total active power output with a moving-average filter and accommodate the losses in the AC model.
The error
\begin{align}\label{eq:err}
e(k)=r-\hat p(k)    
\end{align}
between the reference power output $r$ and the filtered value $\hat p(k)$ is then used as the input for the controller. 
Output of the controller $\pi(k)$ is a function of error $e(k)$ and an inner state of the controller $x_c (k)$ and 
the deadband $\delta$ \eqref{eq:deadband}. 
The control signal is given by the PI regulator or its lag approximant, respectively, as 
\begin{subequations}\label{eq:two-contrlr}
\begin{smalleralign}[\scriptsize]
\pi_{\text{PI}}(k+1)&=
\begin{cases}
K_p e(k)+K_i \left(x_c (k)+e(k)\right), & |e(k)|\geq\delta, \\
\pi_{\text{PI}}(k), & |e(k)|<\delta, 
\end{cases} \label{eq:two-contrlr_pi}
\\
\pi_{\text{Lag}}(k+1)&=
\begin{cases}
K_p e(k)+K_i \left(0.99x_c (k)+e(k)\right), & |e(k)|\geq\delta, \\
\pi_{\text{Lag}}(k), & |e(k)|<\delta.
\end{cases} \label{eq:two-contrlr_lag}
\end{smalleralign}    
\end{subequations}

The procedure is demonstrated by pseudocode in Algorithm \ref{alg:sim}.
\begin{algorithm}[bt]
\SetAlgoLined
$P(0) = \textrm{power-flow}(u(0))$ for a given initial commitment $u(0)$ \\
\For{$i\gets0$ \KwTo $\infty$, i.e., each run of the simulation }{ 
\For{$i\gets1$ \KwTo $k_{\max}$, where the length $k_{\max}$ of the time horizon is given}{ 
		$p(k) = \sum P(k-1)$, cf. \eqref{eq:output-gnrtr} \\
		$\hat p(k) = {\mathcal F}(p(k))$, e.g. \eqref{eq:filtr}  \\
		$e(k) = r - \hat p(k)$, cf. \eqref{eq:err}  \\
            \If{$|e(k)| \geq \delta$}{
                $\pi(k) = {\mathcal C}(e(k))$, cf. \eqref{eq:two-contrlr} \\
            \Else{$\pi(k) = \pi(k-1)$}}
		$u(k) = {\mathcal S}(\pi(k))$, cf. \eqref{eq:prob-func} or \eqref{eq:prob-func2}  \\
		$P(k) = \textrm{power-flow}(u(k))$, cf. \eqref{eq:output-gnrtr}
	}
}
 \caption{Pseudo-code capturing the closed loop.}
 \label{alg:sim}
\end{algorithm}

\begin{figure}[t!]
    \includegraphics[width=\columnwidth]{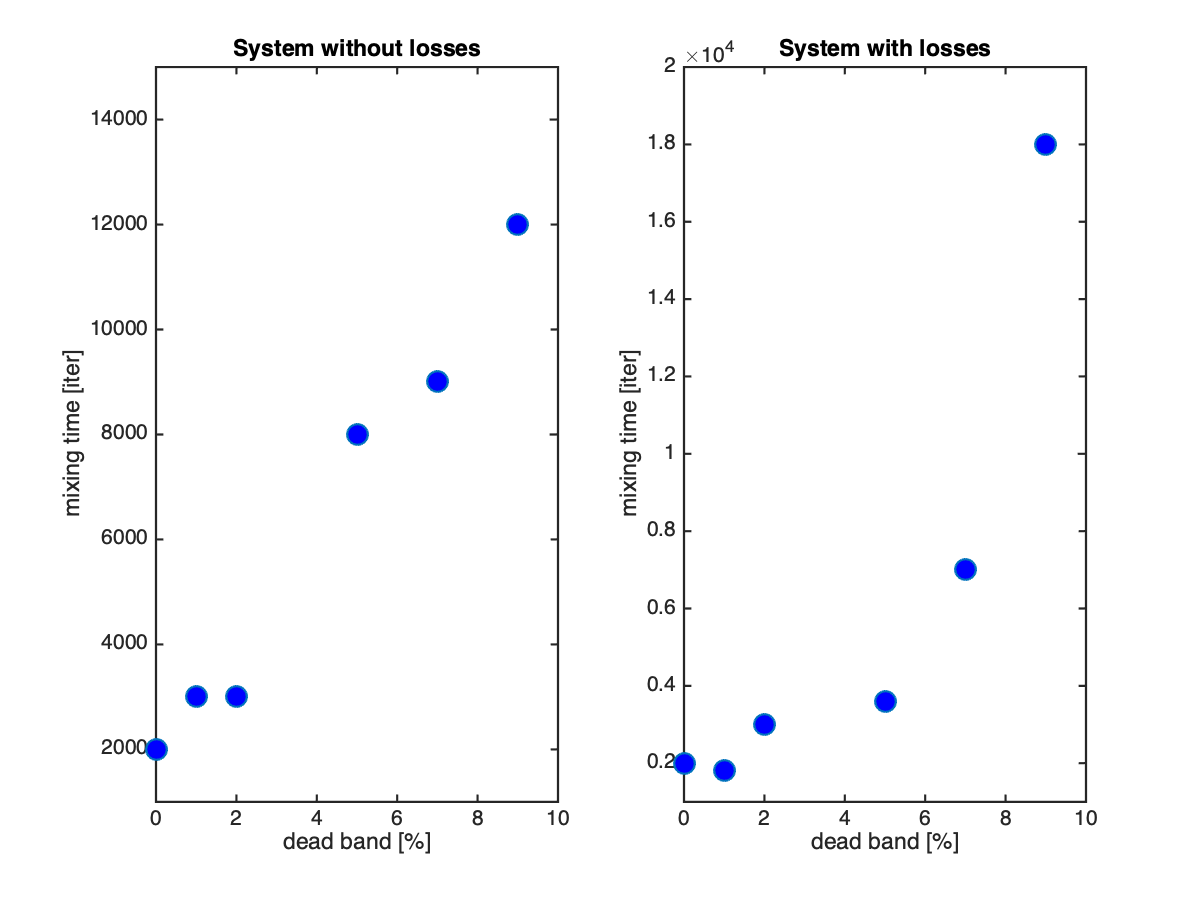}
    \caption{Mixing time as function of deadband expressed in percentage of reference value $r$.}
    \label{fig:mixing_time}
    \end{figure}

For the first test case, we repeated the simulations 50 times, considering various time horizon ranging between 5000 and 80000 time steps. Various sizes of deadband were tested, with deadband being a percentage of reference value $r$. In Figure  \ref{fig:mixing_time}, the size of deadband is shown to affect mixing time of the system, i.e., the number of steps the system needs to reach a certain neighbourhood of the unique invariant measure. In our case, we check if all pairs of initial states $x_{c_A}(k=0), x_{c_B}(k=0)$ of the controller, which we test, are within $0.01\Delta$ with respect to the 2-Wasserstein distance:

\begin{align} \label{eq:Wass}
    W(x_{c_A},x_{c_B}) = \int_0^1 \lvert \hat F_A^{-1}(q) - \hat F_B^{-1}(q) \rvert^2dq \leq 0.01\Delta
 \end{align}

 where $\Delta = W(x_{c_A}(k=1),x_{c_B}(k=1))$ and $\hat F_B, \hat F_B$ are estimates of the distribution functions of $x_{c_A}$ and $x_{c_B}$.

\begin{figure}[t!]
\centering 
No deadband:\\
\includegraphics[width=0.45\columnwidth]{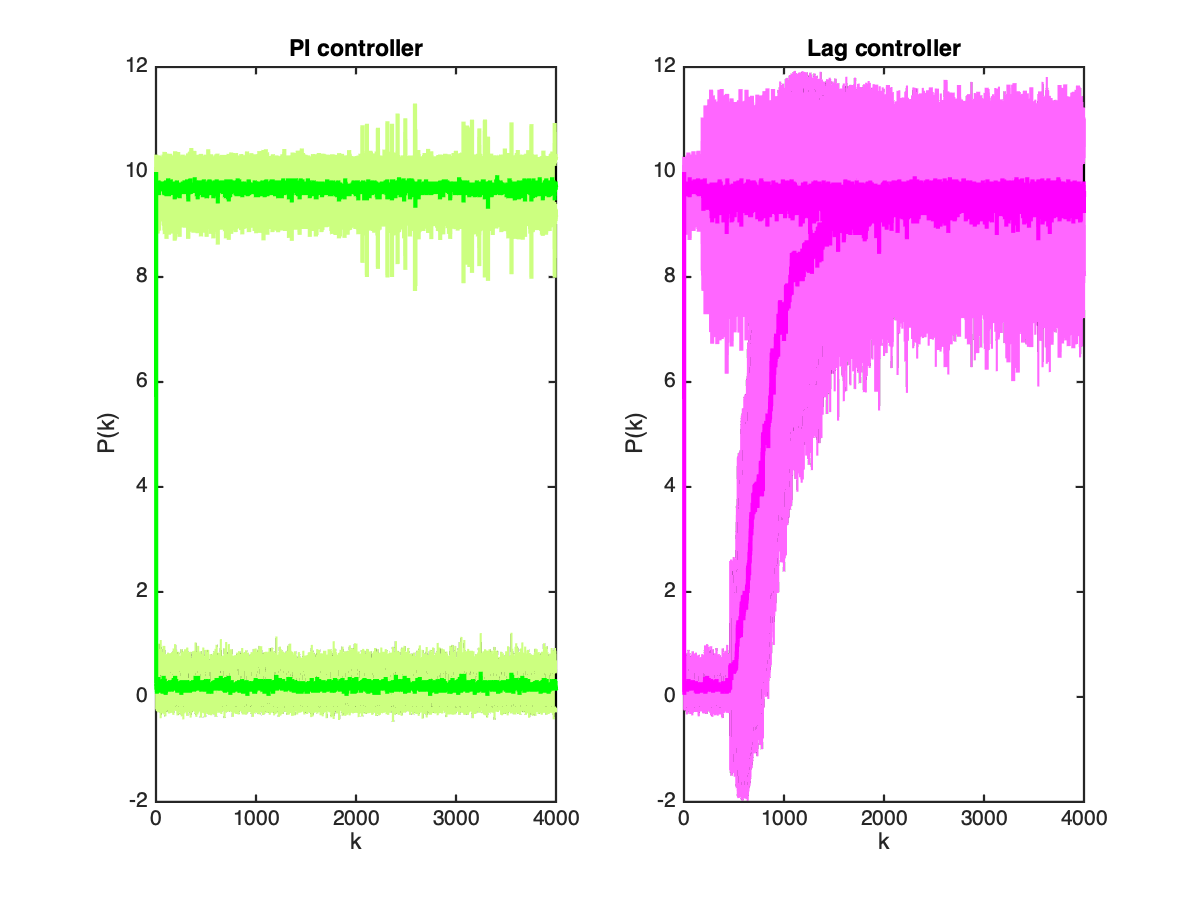}
\includegraphics[width=0.45\columnwidth]{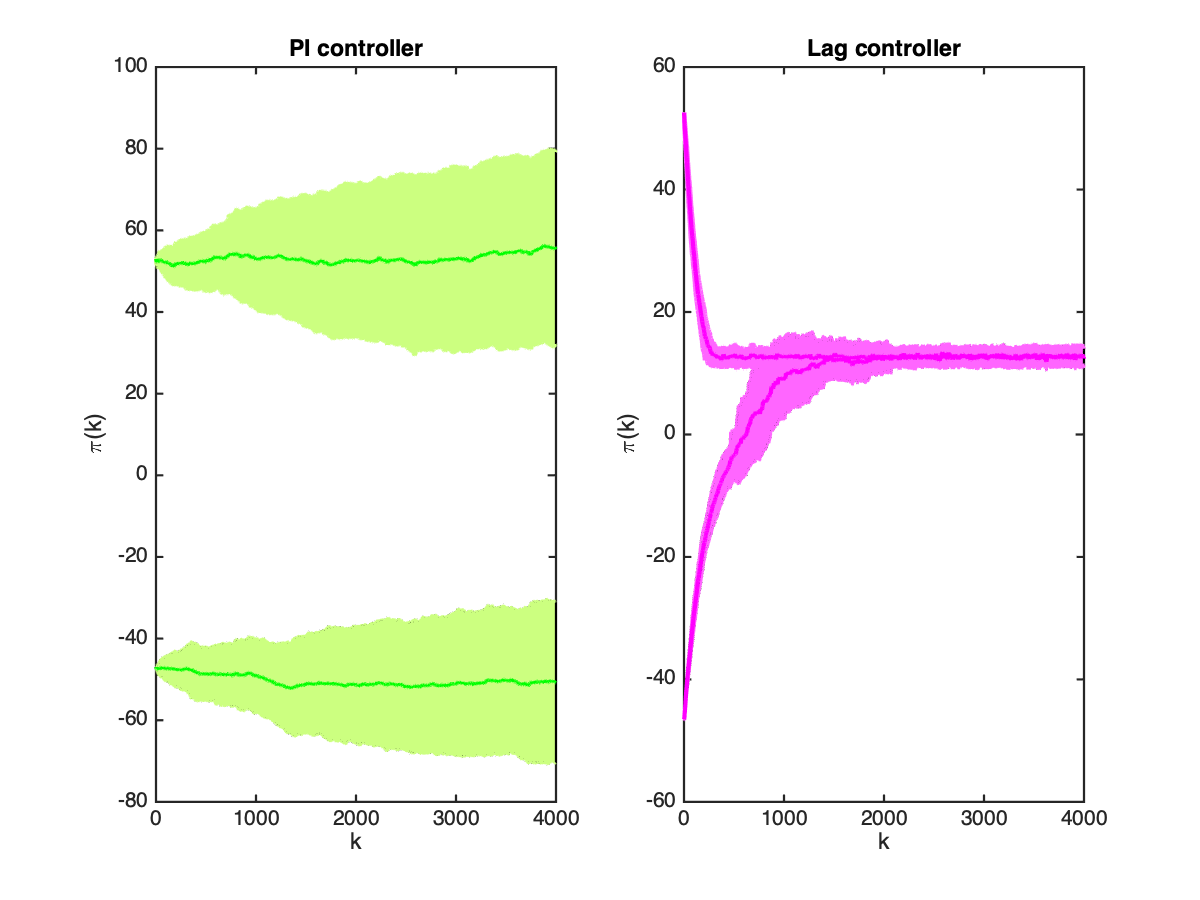}\\
Deadband of 5\% of $r$:\\
\includegraphics[width=0.45\columnwidth]{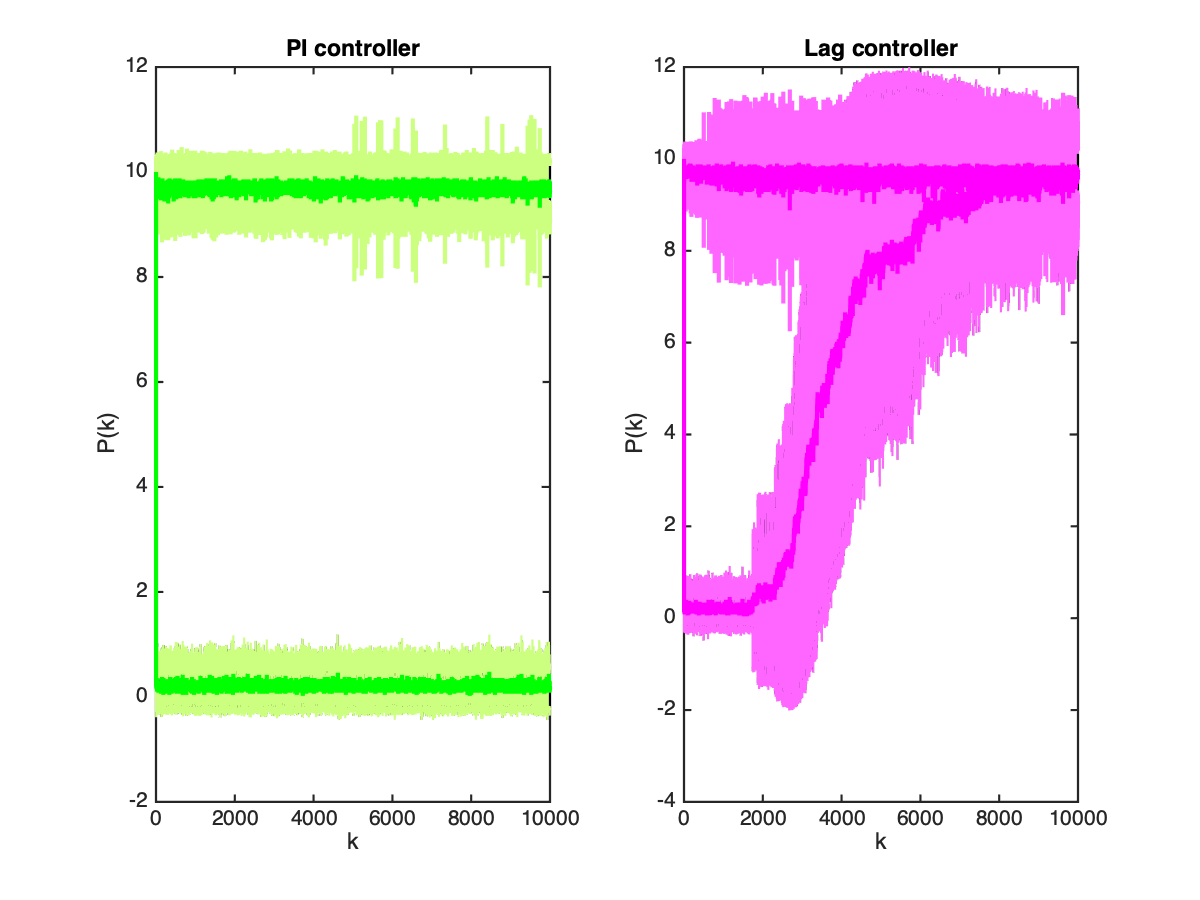}
\includegraphics[width=0.45\columnwidth]{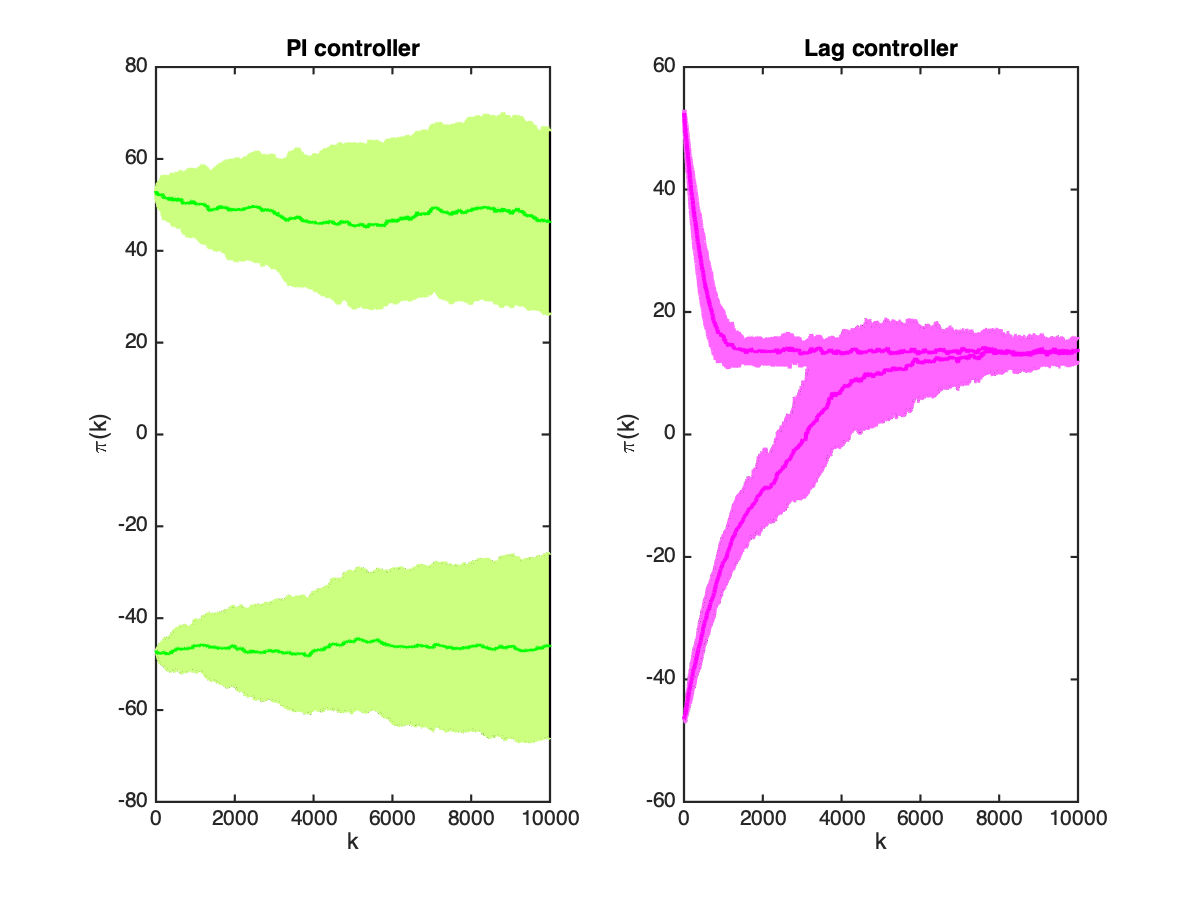}\\

\caption{Results of simulations on the first test case:
The state of the controllers (right) 
and the state of agents (left) utilising probability function $g_{i1}$ of \eqref{eq:prob-func}
as functions of time, for the two controllers and two initial states of each of the two controllers.}
\label{fig:deadband-nolosses}
\end{figure}

\begin{figure}[t!]
\centering 
No deadband:\\
\includegraphics[width=0.45\columnwidth]{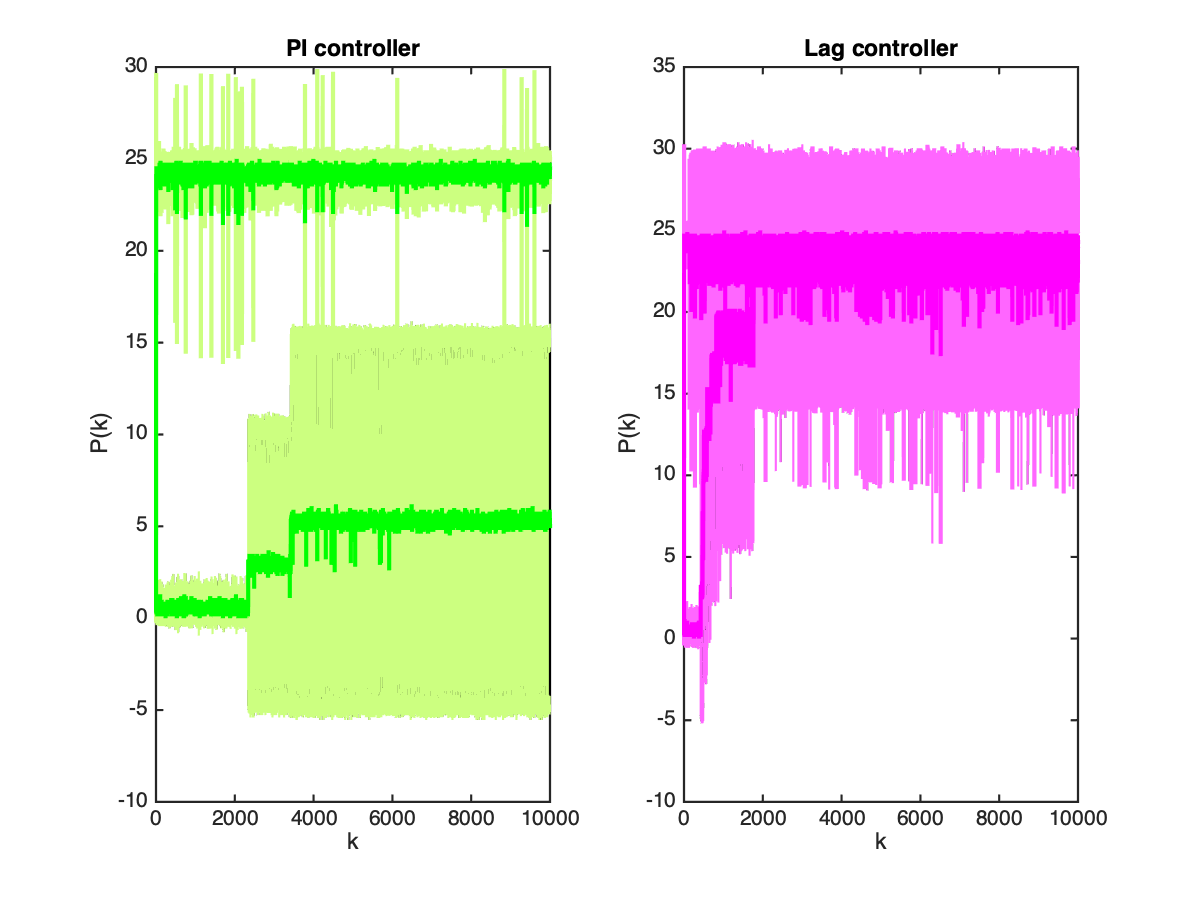}
\includegraphics[width=0.45\columnwidth]{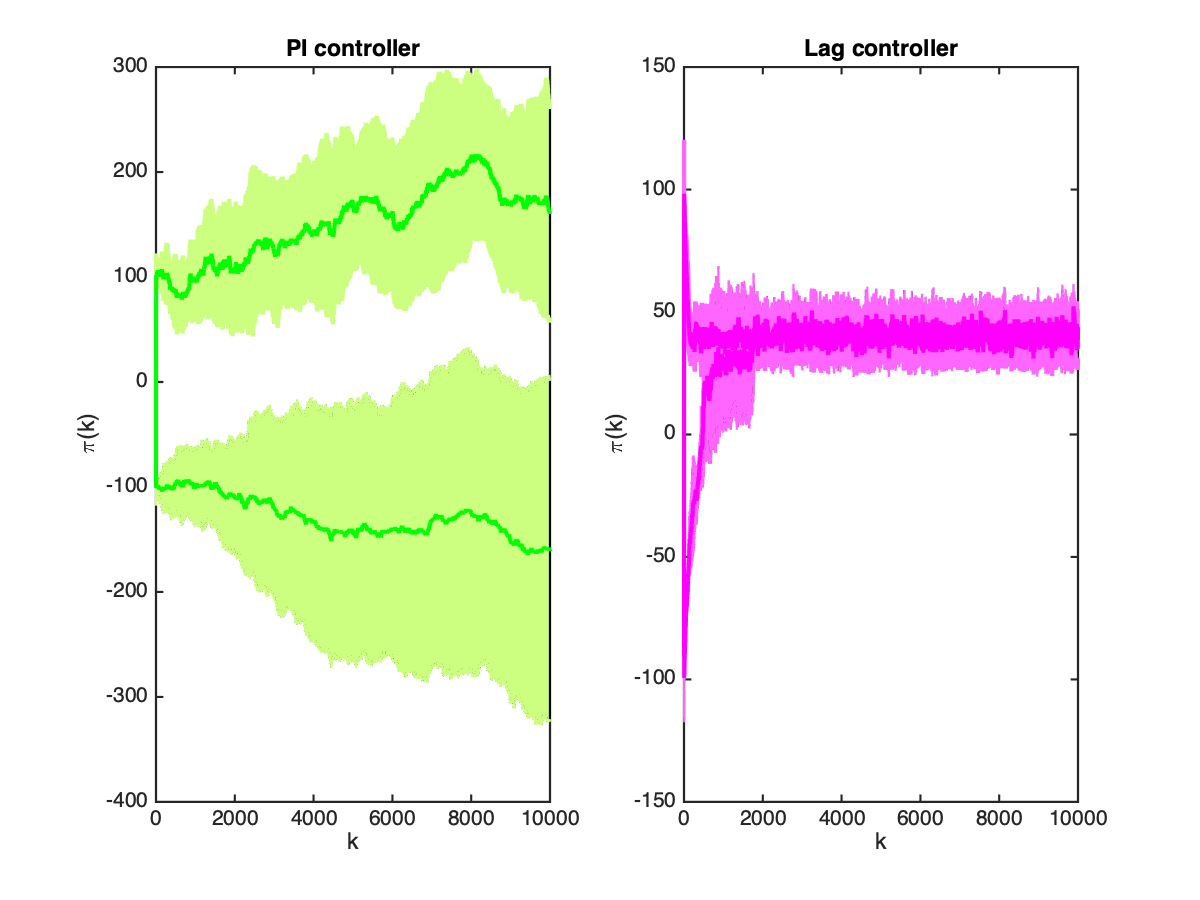}\\
Deadband of 5\% of $r$:\\
\includegraphics[width=0.45\columnwidth]{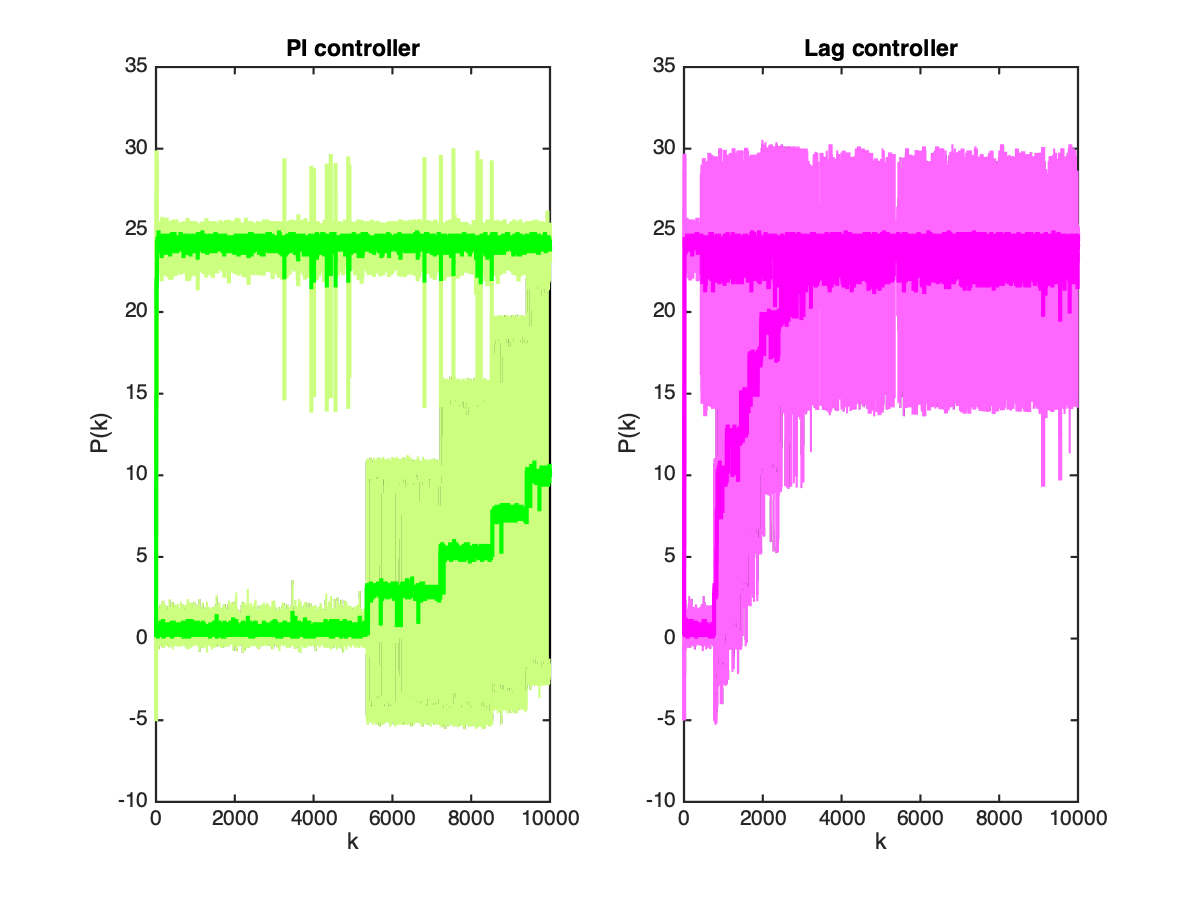}
\includegraphics[width=0.45\columnwidth]{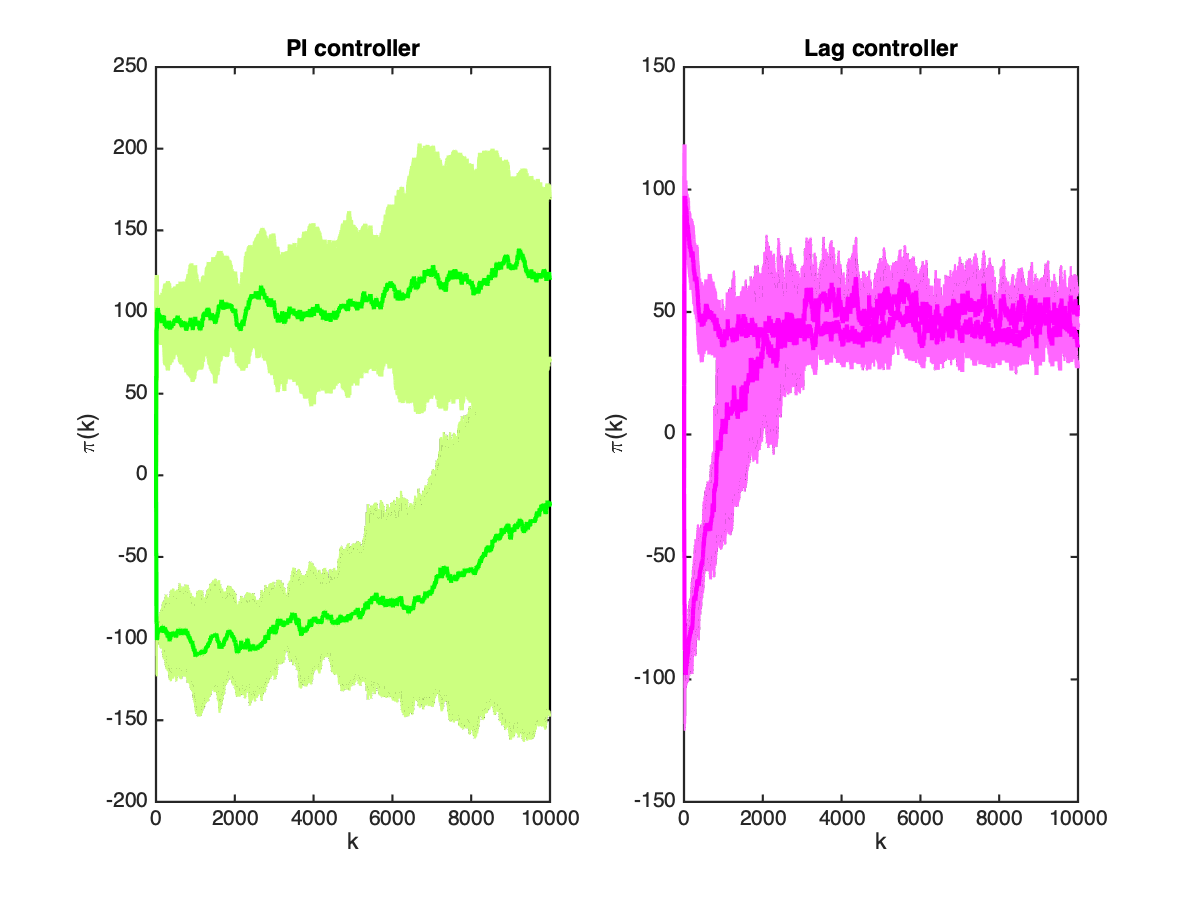}\\

\caption{Results of simulations on Standard IEEE 30-bus test case:
The state of the controllers (right) 
and the state of agents (left) utilising probability function $g_{i1}$ of \eqref{eq:prob-func}
as functions of time, for the two controllers and two initial states of each of the two controllers.}
\label{fig:deadband-losses}
\end{figure}

Still for the second test case, we repeated the simulation 50 times for varying time horizon from 8000 to 100000 time steps. The results in Figure \ref{fig:deadband-nolosses} present the mean over the 50 runs with a solid line and the region of mean $\pm$ standard deviation as a shaded area. Similarly, results for the MATPOWER test case with and without deadband are presented in Figure \ref{fig:deadband-losses}.

In both test cases, systems behave similarly to results in \cite{marecek2021predictability}. Both types of controller (PI and lag approximant) are able to provide regulation to the reference value $r$. Only lag approximant, however, provides ergodic regulation in accordance with Section \ref{ssec:model}.
 Altered behavior of system with implemented discontinuities can be seen in Figure \ref{fig:mixing_time}, where the mixing time increases with the size of deadband, due to the fact that responsiveness of the system is decreased by lack of response within the deadband interval.

\section{Conclusions and Further Work}

We have studied the impact of the use of deadband on the notions of predictability of fairness in load aggregation and demand-response management, which relies on the existence of a unique (Filippov) invariant measure for a stochastic model of a closed-loop model of the system.

\paragraph*{Acknowledgements}
This work received funding from the National Centre for Energy II (TN02000025). FVD acknowledges the support of INdAM-GNCS.

\bibliography{ref,mdps,power,disc}

\FloatBarrier

\end{document}

\clearpage
\onecolumn
\normalsize
\appendix

For the convenience of the reader, we attach an overview of the notation and additional numerical results in the following appendices. 

\subsection{An Overview of Notation}

\begin{tabularx}{\linewidth}{ l | X }
\caption{A Table of Notation}\\\toprule\endfirsthead
\toprule\endhead
\midrule\multicolumn{2}{l}{\itshape continues on next page}\\\midrule\endfoot
\bottomrule\endlastfoot
\textbf{Symbol}                        & \textbf{Meaning}\\\midrule
$\mathbb{N}$                           & the set of natural numbers.\\
$\mathbb{Q}$                           & the set of rational numbers.\\
$\mathbb{R}$                           & the set of real numbers.\\
${\mathbb X}_i$                        & states-space of agent $i$.\\
${\mathbb X}_{\mathcal C}$             & states-space of controller.\\
${\mathbb X}_{\mathcal F}$             & states-space of filter.\\
${\mathbb X}$                          & overall states-space of the closed loop system.\\
${\mathbb Z}$                          & set of all integers.\\
$\mathcal G$                           & an event i.e an element in $\mathcal{B}\left(\Sigma\right)$.\\
$\mathcal J$                           & an index set.\\
$\Sigma$                               & a closed subset of $\mathbb R^n$.\\
$\mathcal{B}\left(\Sigma\right)$       & A Borel $\sigma$ algebra on $\Sigma$, consisting of all possible events.\\
$\{X_k\}_{k\in\N}$                     & A discrete-time Markov chain with state-space $\mathcal K$.\\
$\mathbf{P}$                           & A transition-probability operator.\\
$\mathbb P_{\lambda}$                  & the probability measure induced on
the path space.\\
$\lambda_k$                            & probability measures.\\
$\mu$                                  & a probability measure.\\
$j$                                    & an element in $\mathcal J$.\\
$\ell$                                 & cardinality of output maps.\\
$f_j$                                  & a function on $\mathcal K$.\\
$p_j$                                  & probability functions.\\
$\pi(k)$                               & control signal.\\
$\Pi$                                  & space of all control signal.\\
$\mathbb A_i$                          & space of all agent's actions.\\
$a_1,a_2,\dots, a_{L_i}$               & possibles actions.\\
$L_i$                                  & cardinality of action space of agent $i$.\\
$\mathbb D_i$                          & demand space of agent $i$.\\
$d_{i,1},\dots, d_{i,m_i}$             & possible demands.\\
$m_i$                                  & cardinality of demand space of agent $i$.\\
$n_i$                                  & dimension of the state-space of $i^{th}$ agent's private state-space.\\
$w_i$                                  & an affine map.\\
$h_i$                                  & number of output maps.\\
$\mathcal H_{i\ell}$                   & output maps.\\
$y_i(k)$                               & output of agent $i$.\\
$\hat y_(k)$                           & value of $y_(k)$ filtered by filter.\\
$p_{i\ell}(\pi)$                       & probability functions.\\
$p'_{i\ell}(\pi)$                      & probability functions.\\
$x_i(k)$                               & internal state of agent $i$ at time $k$.\\
$\mathcal W_{ij}$                      & transition maps.\\
$\mathcal F$                           & a filter.\\
$\mathcal C$                           & a controller.\\
$x_f(k)$                               & the internal state of the filter at time $k$ of dimension $n_f$.\\
$x_c(k)$                               & the internal state of the filter at time $k$ of dimension $n_c$\\
$\mathcal W_{c}$                       & map moddeling the controller.\\
$\mathcal H_{c}$                       & map moddeling the controller\\
$\epsilon$                             & a sufiiciently small positive real number\\
$l_m$                                  & Lipschitz constants.\\
$\nu^{\star}$                          & an invariant distribution.\\
$\beta$                                & a class $\mathcal {K L}$ function.\\
$\gamma$                               & a class $\mathcal K$ function.\\
$u(k)$                                 & a binary vector indicating the commitment of DERs in \eqref{eq:cmmtd-gnrtr}.\\
$p(k)$                                 & filtered value of the active power output of DERs in \eqref{eq:filtr}.\\
$P(k)$                                 & active power output of DERs in \eqref{eq:output-gnrtr}.\\
$\hat p(k)$                           &  aggregate active power output of the ensemble at time $k$.\\
$e(k)$                                 & error defined in \eqref{eq:err}.
\end{tabularx}

\end{document}